\documentclass[12pt]{amsart}
\usepackage{amssymb}
\usepackage[all]{xy}
\UseComputerModernTips
\newtheorem{theorem}{Theorem}[section]
\newtheorem{lemma}[theorem]{Lemma}
\newtheorem{proposition}[theorem]{Proposition}
\newtheorem{corollary}[theorem]{Corollary}
\theoremstyle{definition}
\newtheorem{definition}[theorem]{Definition}

\newtheorem{example}[theorem]{Example}

\theoremstyle{remark}
\newtheorem{remark}[theorem]{Remark}

\newcommand{\C}{{\mathcal C}}
\newcommand{\D}{{\mathcal D}}
\newcommand{\E}{{\mathcal E}}
\newcommand{\Oh}{{\mathcal O}}
\newcommand{\bndry}{\partial}
\DeclareMathSymbol{\boxprod}{\mathbin}{AMSa}{"03} 
\newcommand{\disjunion}{\sqcup}
\DeclareMathOperator{\Hom}{Hom}
\newcommand{\iso}{\cong}
\DeclareMathOperator{\Map}{Map}
\newcommand{\Over}{{\downarrow}}
\newcommand{\union}{\cup}
\newcommand{\Union}{\bigcup}
\newcommand{\Cats}{{\mathrm{Cat}}}
\newcommand{\SSets}{{\mathcal S}}
\newcommand{\horn}[2]{\Lambda^{#1}[{#2}]}

\newcommand{\Dash}[2]{\text{$#1$-$#2$}}
\newcommand{\OCats}{\Dash{\Oh}{\Cats}}
\newcommand{\OSSets}{\Dash{\Oh}{\SSets}}
\DeclareMathOperator{\Sd}{Sd}
\DeclareMathOperator{\Ex}{Ex}
\begin{document}
%
\title{Simplicial and categorical comma categories}

\author{Steven R. Costenoble}
\address{Department of Mathematics\\103 Hofstra University\\
   Hempstead, NY 11549-1030}
\email{Steven.R.Costenoble@Hofstra.edu}


\date{November 29, 2000}

\abstract
We consider four categories: the category of diagrams of small categories
indexed by a given small category $\Oh$,
the (comma) category of small categories over $\Oh$,
the category of diagrams of simplicial sets indexed by $\Oh$,
and the category of simplicial sets over the nerve of $\Oh$.
Fritsch and Golasi\'nski claimed that these four categories have
equivalent homotopy categories but, in fact, their proof contains
an error and the homotopy categories are not equivalent with
the weak equivalences they use in the comma categories.
We show here that
the correct weak equivalences are the ``weak fibre homotopy equivalences''
defined by Latch.
We also construct a model category structure on the category
of simplicial sets over $N\Oh$ in which the weak equivalences
are the weak fibre homotopy equivalences.
\endabstract

\maketitle
%

\section{Introduction}

Let $\Cats$ denote the category of small categories,
let $\SSets$ denote the category of simplicial sets, and
let $\Oh$ be a small category.
In \cite{FG:diagrams}, Fritsch and Golasi\'nski considered the
relationship between four categories:
$\OCats$, the category of contravariant functors $\Oh\to\Cats$,
$\Cats\Over\Oh$, the category of small categories over $\Oh$,
$\OSSets$, the category of contravariant functors $\Oh\to\SSets$, and
$\SSets\Over N\Oh$, the category of simplicial sets over the nerve $N\Oh$.
These categories fit into the following diagram.
\[
 \xymatrix{
  \Cats\Over\Oh \ar@<-.25ex>@_{->}[r]_F \ar[d]_N &
  \OCats \ar@<-.25ex>@_{->}[l]_E \ar[d]^N \\
  \SSets\Over N\Oh \ar@<-.25ex>@_{->}[r]_F &
  \OSSets \ar@<-.25ex>@_{->}[l]_E
  }
\]
In this diagram, the maps labeled $F$ are left adjoints of the maps labeled $E$,
the maps labeled $N$ are induced in the obvious way by the nerve functor, and
$FN = NF$.

The categories $\OSSets$ and $\OCats$ arise naturally in equivariant algebraic topology
as models for the homotopy theory of $G$-spaces when
$G$ is a discrete group and $\Oh$ is the category of orbits of $G$.
In that context, the appropriate model category structure on
$\OSSets$ is the one in which a natural transformation
$\phi\to\psi$ is a fibration or a weak equivalence exactly when it is
an objectwise fibration or weak equivalence, respectively.
Similarly, the appropriate weak equivalences in $\OCats$ are the objectwise weak
equivalences (where weak equivalence of small categories is weak equivalence
of their nerves); we can define a model category structure on $\OCats$
in which the fibrations are objectwise fibrations, using the model
category structure given by Thomason \cite{Tho:cat}.
Fritsch and Golasi\'nski claim that the adjunctions above are Quillen
equivalences when we use on $\Cats\Over\Oh$ and $\SSets\Over N\Oh$ 
the model category structures
in which a map over $\Oh$ or $N\Oh$ is a cofibration, fibration, or weak
equivalence exactly when it is so considered as an ordinary
map. Unfortunately, their proof contains an error,
and their claim is not true, as shown by simple examples
(see Example~\ref{ex:error}).

The correct weak equivalences to use in $\Cats\Over\Oh$ and
$\SSets\Over N\Oh$ are the {\em weak fibre homotopy equivalences}
defined in \cite{La:fibred}.
In Section~\ref{sec:equiv} we give a correct proof of the equivalences
of the homotopy categories of $\Cats\Over\Oh$ and $\OCats$ and
of $\SSets\Over N\Oh$ and $\OSSets$. This gives another proof of the equivalence
of the homotopy categories of $\Cats\Over\Oh$ and $\SSets\Over N\Oh$
shown in \cite{La:fibred}.
In the remainder of the paper we construct a model category structure
on $\SSets\Over N\Oh$ having weak fibre homotopy equivalences as its
weak equivalences.
Defining a model category structure on $\Cats\Over\Oh$
appears to be much more difficult and we shall not attempt it here.

\section{The adjoint functors}\label{sec:functors}

In this section we define and examine the
adjoint pairs $(F,E)\colon \Cats\Over \Oh \to \OCats$ and
$(F,E)\colon \SSets\Over N\Oh \to \OSSets$. 

First, let us recall some standard notations.
Let $\Delta$ be the category whose objects are the
ordered sets
$[n] = \{0, 1, \dots, n\}$, $n\geq 0$, and whose maps are the order-preserving
functions. We often think of $[n]$ as a category in the usual way,
so that the maps in $\Delta$ can be thought of as functors.
We write $\Delta[n]$ for the nerve of $[n]$, the standard $n$-simplex.
We write $\bndry\Delta[n]\subset\Delta[n]$ for the subsimplicial
set generated by all of the $(n-1)$-dimensional faces of $\Delta[n]$.
If $0\leq k\leq n$, we write
$\horn kn \subset \Delta[n]$ for the $k$th {\em horn},
the subsimplicial set generated by all of the $(n-1)$-dimensional
faces of $\Delta[n]$ except the one opposite the $k$th vertex.

Now, we define $F\colon \Cats\Over\Oh \to \OCats$ as follows.
If $\phi\colon \C\to\Oh$ is a small category over $\Oh$ and
$b$ is an object of $\Oh$, we let $(F\phi)(b)$ be the pullback in the
following diagram.
\[
 \xymatrix{
  (F\phi)(b) \ar[r] \ar[d] & \C \ar[d]^\phi \\
  b\Over\Oh \ar[r] & \Oh
  }
\]
Here, $b\Over\Oh$ is the usual comma category of objects under $b$.
(The category $(F\phi)(b)$ is otherwise known as $b\Over\phi$.)
If $\beta\colon b\to c$ is an arrow in $\Oh$, we have a functor
$c\Over\Oh\to b\Over\Oh$ given by composition with $\beta$, and this
makes $F\phi$ into a contravariant functor.

The right adjoint $E\colon \OCats\to\Cats\Over\Oh$ is given by
the Grothendieck construction \cite{Gr:groupefond}.
If $\psi\colon\Oh\to\Cats$ is a contravariant functor, we define
$E\psi\colon \bar E\psi\to\Oh$ as follows. The category
$\bar E\psi$ has as set of objects the disjoint union of the sets of objects
of all of the categories $\psi(b)$; $E\psi$ maps each object in $\psi(b)$ to $b$.
If $x\in\psi(b)$ and $y\in\psi(c)$,
a map in $\bar E\psi$ from $x\to y$ is a pair $(\beta,f)$ where
$\beta\colon b\to c$ in $\Oh$ and $f\colon x\to\beta^* y$
in $\psi(b)$. Composition of pairs is given by
$(\beta,f)\circ(\gamma,g) = (\beta\gamma, \gamma^* f\circ g)$.

We can check directly that $F$ and $E$ are adjoint, or we can do the following.
Note first that, if $\phi\colon \C\to\Oh$, then
\[
 \phi \iso \int^{\sigma\in\Delta\Over\Oh} \!\!\!\!\Hom_\Oh(\sigma,\phi)\times\sigma
\]
where $\Delta\Over\Oh$ is the full subcategory of $\Cats\Over\Oh$ containing
all maps $[n]\to\Oh$, and $\Hom_\Oh$ denotes the set of maps over $\Oh$.
(In fact, we could restrict to the subcategory of $\Delta\Over\Oh$ containing
only those objects $[n]\to\Oh$ with $n\leq 2$.)
It is easy to see that $F$ preserves colimits, so we have that
\[
 F\phi \iso \int^\sigma \!\!\Hom_\Oh(\sigma,\phi)\times F\sigma.
\]
It follows from general nonsense about coends that $F$ then has as its
right adjoint the functor $E$ defined by
\[
 E\psi = \int^\sigma \!\!\Hom_\Oh(F\sigma,\psi)\times \sigma,
\]
where $\Hom_\Oh$ here denotes the set of natural transformations of
functors defined on $\Oh$. Remembering that we need only consider
$n\leq 2$, it's a pleasant exercise to see that $E$ so defined can
be described by the Grothendieck construction.

The constructions for simplicial sets are similar.
The functor $F\colon \SSets\Over N\Oh\to \OSSets$
is defined as follows. Let $\phi\colon X\to N\Oh$ be a map of simplicial
sets and let $b$ be an object of $\Oh$. We let $(F\phi)(b)$ be the pullback in
the following diagram.
\[
 \xymatrix{
  (F\phi)(b) \ar[r] \ar[d] & X \ar[d]^{\phi} \\
  N(b\Over\Oh) \ar[r] & N\Oh
  }
\]
We make $(F\phi)(b)$ a contravariant functor of $b$ in the evident way.
(Note: Since $N$ preserves pullbacks, this diagram show that $FN=NF$.)

To define the right adjoint $E$, we first notice that $F$ preserves
colimits, as may be checked directly. Now, let $\phi\colon X\to N\Oh$ and
$\sigma\colon\Delta[n]\to N\Oh$.
(We shall freely move back and forth between such a map $\sigma$ and
the corresponding map $[n]\to\Oh$, which we shall also call $\sigma$.)
Let $X_\sigma\subset X_n$ be the
set of $n$-cells of $X$ that map to $N\Oh$ via $\sigma$.
Put another way, $X_\sigma = \Hom_{N\Oh}(\sigma,\phi)$.
We then have
\begin{equation*}
 \phi \iso \int^{\sigma\in\Delta\Over \Oh} \!\!\!\!X_\sigma \times \sigma.
\end{equation*}
(In fact, it is often convenient to think of a simplicial set over $N\Oh$
as instead a contravariant set-valued functor on $\Delta\Over \Oh$.)
Using the fact that $F$ preserves colimits, we can then write
\begin{equation*}
 F\phi \iso \int^\sigma \!\!X_\sigma\times F\sigma.
\end{equation*}
It follows that the right adjoint $E$ is given as follows.
If $\psi\colon \Oh\to \SSets$ is a contravariant functor,
then
\[
 E\psi = \int^{\sigma\in\Delta\Over\Oh} \!\!\!\!\Hom_\Oh(F\sigma,\psi)\times\sigma,
\]
so that
\begin{equation*}
 (E\psi)_\sigma = \Hom_{\Oh}(F\sigma,\psi)
\end{equation*}
and
\begin{equation*}
 (E\psi)_n = \Union_{\sigma\colon\Delta[n]\to N\Oh}
   \!\!\!\!\Hom_{\Oh}(F\sigma,\psi).
\end{equation*}
This is a generalization to
simplicial sets of the Grothendieck construction for categories.
As noted in \cite{FG:diagrams}, it is also the homotopy colimit of $\psi$.

\begin{example}\label{ex:error}
(See also Remark~3.10 in \cite{La:fibred}.)
Contrary to the claim in \cite{FG:diagrams},
if $\phi\colon \C\to \Oh$,
$\psi\colon \D\to\Oh$, and $f\colon\phi\to\psi$ is a functor such that
$f\colon\C\to\D$ is a weak equivalence of categories, then $Ff$ need not be
an objectwise weak equivalence.
Consider the following example. Let
$\Oh = [1]$, let $\C = [0]$ with $\phi\colon\C\to\Oh$ taking 0 to 0,
and let $\D=[1]$ with $\psi\colon\D\to[1]$ the identity.
Let $f\colon\phi\to\psi$ take 0 to 0.
Then $f$ is a weak equivalence of categories but
$(F\phi)(1)$ is empty and $(F\psi)(1)$ is not, so $Ff$ cannot
be an objectwise weak equivalence.
Taking the nerve gives an example showing that the similar claim
about simplicial sets is also false.
In fact, we shall see below that this example is prototypical of the
cases that we need to avoid.
\end{example}

We now examine the fundamental
diagrams $F\sigma$ in detail in the simplicial case.
First, notice that $\sigma\colon [n]\to \Oh$ is determined by
a sequence of maps
\begin{equation*}
 \sigma_0 \xrightarrow{g_1} \sigma_1 \xrightarrow{g_2} 
   \cdots \xrightarrow{g_n} \sigma_n
\end{equation*}
in $\Oh$. 
If $b$ is an object of $\Oh$, we have the following pullback diagram.
\[
 \xymatrix{
  (F\sigma)(b) \ar[r] \ar[d] & \Delta[n] \ar[d]^\sigma \\
  N(b\Over\Oh) \ar[r] & N\Oh
  }
\]
Thus, a $k$ simplex in $(F\sigma)(b)$ can be described as a pair $(\alpha,\omega)$
where $\alpha\colon [k]\to[n]$ and $\omega\colon [k]\to b\Over\Oh$,
with
\begin{equation*}
 \omega = (b\to \sigma\alpha(0) \xrightarrow{\sigma\alpha(0\to 1)}
   \sigma\alpha(1)\xrightarrow{\sigma\alpha(1\to 2)}\cdots
   \xrightarrow{\sigma\alpha(k-1\to k)}\sigma\alpha(k)).
\end{equation*}
Write $\Oh(b,c)$ for the set of maps from $b$ to $c$
in $\Oh$ and, for a fixed object $c$,
write $\Oh(-,c)$ for the evident contravariant set-valued functor
on $\Oh$.
Define $\iota_n\colon \Oh(-,\sigma_0)\times\Delta[n]\to F\sigma$ to be the map taking
$(g,\alpha)\in \Oh(b,\sigma_0)\times\Delta[n]_k$ to $(\alpha,\omega)$, where
\begin{equation*}
 \omega = (b\xrightarrow{\sigma(0\to\alpha(0))\circ g}
          \sigma\alpha(0)\xrightarrow{\sigma\alpha(0\to 1)}\sigma\alpha(1)
          \xrightarrow{\sigma\alpha(1\to 2)}\cdots
          \xrightarrow{\sigma\alpha(k-1\to k)} \sigma\alpha(k)).
\end{equation*}

\begin{lemma}\label{lem:attachment}
Let $\sigma\colon\Delta[n]\to N\Oh$, let
$A\subset\Delta[n]$ be a subsimplicial set containing
the face opposite $0$, and let $\lambda$ denote the composite
$A\to\Delta[n]\to N\Oh$. Then the following is a pushout diagram,
where the top map is the restriction of $\iota_n$.
\[
 \xymatrix{
  \Oh(-,\sigma_0)\times A \ar[r] \ar[d] & F\lambda \ar[d] \\
  \Oh(-,\sigma_0)\times \Delta[n] \ar[r]_-{\iota_n} & F\sigma
  }
\]
\end{lemma}

\begin{proof}
Let $(\alpha,\omega)$ be an element of $(F\sigma)(b)$ for some object $b$.
Then, $(\alpha,\omega)\in F\lambda$ if and only if $\alpha\in A$.
Let us assume that $\alpha\notin A$. Since $A$ contains the face opposite 0,
we must have that $\alpha(0) = 0$.
Then $\omega$ has the form
\begin{equation*}
 \omega = (b\xrightarrow{g}\sigma_0\to c_1\to\cdots\to c_k),
\end{equation*}
and so $(\alpha,\omega)$ is the image under $\iota_n$ of a unique element of
$\Oh(b,\sigma_0)\times\Delta[n]$, i.e., $(g,\alpha)$.
Thus, $F\sigma$ is the pushout claimed.
\end{proof}

We can now describe the structure of $F\sigma$.

\begin{proposition}\label{prop:structure}
If $\sigma\colon\Delta[0]\to N\Oh$, then
$\iota_0\colon \Oh(-,\sigma_0)\times\Delta[0] \to F\sigma$
is an isomorphism.
If $\sigma\colon\Delta[n]\to N\Oh$ with $n>0$, we have the following
pushout diagram.
\[
 \xymatrix{
  \Oh(-,\sigma_0)\times\Delta[n-1] \ar[r]^-{\iota_{n-1}g_{1*}} \ar[d]_{\delta_0}
    & F(\sigma\delta_0) \ar[d] \\
  \Oh(-,\sigma_0)\times\Delta[n] \ar[r]_-{\iota_n} & F\sigma
  }
\]
Here, $\delta_0\colon\Delta[n-1]\to\Delta[n]$ is inclusion of the face
opposite $0$.
\end{proposition}

\begin{proof}
This follows from Lemma~\ref{lem:attachment} on letting
$A$ be the face opposite 0 ($A=\emptyset$ if $n=0$).
\end{proof}

This gives us a nice, Grothendieck-like interpretation of the right adjoint
$E$. If $\psi\colon\Oh\to\SSets$ and
$\sigma\colon\Delta[n]\to N\Oh$,
the description of $F\sigma$ as a pushout allows us to describe
$(E\psi)_\sigma = \Hom_{\OSSets}(F\sigma,\psi)$ as a pullback. 
Interpreting the terms involved, we get the following.

\begin{corollary}\label{cor:fsigma}
Let $\psi\colon\Oh\to\SSets$ and $\sigma\colon\Delta[n]\to N\Oh$,
with
\begin{equation*}
 \sigma = (\sigma_0\xrightarrow{g_1}\cdots\xrightarrow{g_n}\sigma_n).
\end{equation*}
Then an element of $(E\psi)_\sigma$ is specified by a sequence
$(\tau_0,\tau_1,\dots,\tau_n)$ where
$\tau_k\in \psi(\sigma_k)_{n-k}$ and
$d_0\tau_{k-1} = g_k^*\tau_k$ for $1\leq k \leq n$.
\end{corollary}

We can elaborate Example~\ref{ex:error} to show that, 
if $\sigma\colon\Delta[n]\to N\Oh$
and $\lambda$ is the composite $\horn 0 n\to\Delta[n]\to N\Oh$, then
$(F\lambda)(b)\to (F\sigma)(b)$ need not be a weak equivalence for every $b$.
However, this is a problem for the 0th horn only.

\begin{proposition}\label{prop:we}
Let $\sigma\colon\Delta[n]\to N\Oh$ and let $\lambda$
be the composite $\horn k n\to\Delta[n]\to N\Oh$ where $k > 0$. Then
$(F\lambda)(b)\to (F\sigma)(b)$ is a weak equivalence for every
object $b$ in $\Oh$.
\end{proposition}

\begin{proof}
Letting $A = \horn k n$ in Lemma~\ref{lem:attachment}
(notice that this requires $k>0$),
we have the following pushout diagram.
\[
 \xymatrix{
  \Oh(b,\sigma_0)\times\horn k n \ar[r] \ar[d] & (F\lambda)(b) \ar[d] \\
  \Oh(b,\sigma_0)\times\Delta[n] \ar[r] & (F\sigma)(b)
  }
\]
But then standard results about simplicial sets say that
$(F\lambda)(b)\to (F\sigma)(b)$ is a weak equivalence.
\end{proof}

Another obvious application of Lemma~\ref{lem:attachment} is to the case where
$A = \bndry\Delta[n]$. We shall use this later
to examine the effect on cell structures of
applying $F$.

\section{The equivalences of homotopy categories}\label{sec:equiv}

Although not usually stated explicitly,
the following result or some variation on it
is often used to show equivalences of homotopy categories when one category
lacks a full model structure.
Recall that, if $\C$ is a category with certain of its maps designated as
weak equivalences, then its homotopy category $h\C$ is the category obtained by
inverting the weak equivalences; this category may or may not exist, the
difficulty being a set-theoretic one
we know is overcome for model categories.

\begin{proposition}\label{prop:equivalence}
Let $\C$ be a category and let $\D$ be a model category.
Let $(L,R)\colon\C\to\D$ be an adjoint pair such that
the counit $\epsilon\colon LR(d)\to d$ is a weak equivalence for every
object $d$ in $\D$.
Say that a map $f\colon x\to y$ in $\C$ is a weak equivalence
when $Lf$ is one.
Then the homotopy category of $\C$ exists and is equivalent to the
homotopy category of $\D$.
\end{proposition}

The dual statement is the following.

\begin{proposition}\label{prop:dualequivalence}
Let $\C$ be a model category and let $\D$ be a category.
Let $(L,R)\colon\C\to\D$ be an adjoint pair such that the
unit $\eta\colon c\to RL(c)$ is a weak equivalence for every
object $c$ in $\C$.
Say that a map $f\colon x\to y$ in $\D$ is a weak equivalence
when $Rf$ is one.
Then the homotopy category of $\D$ exists and is equivalent to the
homotopy category of $\C$.
\end{proposition}

The proofs of these results are straightforward. We can weaken
the hypotheses: In Proposition~\ref{prop:equivalence} it suffices to
assume that $\epsilon$ is a weak equivalence only on fibrant objects,
and in Proposition~\ref{prop:dualequivalence} it suffices to assume
that $\eta$ is a weak equivalence only on cofibrant objects.
Of course, if both categories are model categories then we can weaken
the hypotheses further to the requirement that $(L,R)$ be a Quillen equivalence.

In the context of Proposition~\ref{prop:dualequivalence},
we can often push the conclusion further: If $\C$ is cofibrantly generated
and certain smallness assumptions hold, then $L$ carries the generating
cofibrations and acyclic cofibrations to cofibrant generators for
a model category structure on $\D$. In that case $(L,R)$ becomes a
Quillen equivalence. Proposition~\ref{prop:dualequivalence}
underlies the various proofs that the homotopy category
of simplicial sets is equivalent to the homotopy category of small
categories.
Illusie~\cite{Illu:complex}, Lee~\cite{Lee:homotopy}, Latch~\cite{La:homology}
and Thomason~\cite{Tho:cat} constructed various adjunctions
$(L,R)\colon\SSets\to\Cats$ whose units are weak equivalences,
with $R$ naturally weakly equivalent to $N$.
Any of these can be used to define a
cofibrantly generated model category structure on
$\Cats$; most of the work in \cite{Tho:cat} (as corrected in
\cite{Cis:dwyer}) is to show that the structure constructed there
is proper.

As mentioned in the introduction, it's well-known and not difficult to show
that the category $\OSSets$ has a model category structure in which the weak
equivalences are the objectwise weak equivalences and the fibrations are
the objectwise fibrations. In fact, this structure is cofibrantly generated:
For a set of generating cofibrations we may take the collection of maps
$\Oh(-,c)\times\bndry\Delta[n]\to\Oh(-,c)\times\Delta[n]$,
and for a set of generating acyclic cofibrations we may take the collection
of maps
$\Oh(-,c)\times\horn kn \to\Oh(-,c)\times\Delta[n]$, for all $k$ and $n$.

Similarly, using the
model category structure on $\Cats$ constructed by Thomason~\cite{Tho:cat},
the category $\OCats$ has a cofibrantly generated
model category structure in which the
weak equivalences are the objectwise weak equivalences and the fibrations
are the objectwise fibrations.
The adjunction $(c\Sd^2,\Ex^2N)\colon \SSets\to\Cats$ is a Quillen equivalence
and is easily seen to induce a Quillen equivalence between
$\OSSets$ and $\OCats$.
Moreover, because $N\to \Ex^2 N$ is a natural weak equivalence, we
have the following result.

\begin{theorem}\label{thm:diagramequivalence}
With weak equivalences defined objectwise, $N\colon\OCats\to\OSSets$
induces an equivalence of homotopy categories.
\end{theorem}

We now define weak equivalence in $\Cats\Over\Oh$ by saying that
a map $f$ over $\Oh$ is a weak equivalence when $Ff$ is one in $\OCats$.
Similarly, we say that a map $f$ in $\SSets\Over N\Oh$ is a weak equivalence
when $Ff$ is one in $\OSSets$. These are the notions of
weak fibre homotopy equivalence from \cite{La:fibred}.

To use Proposition~\ref{prop:equivalence}, we need to show that
the counits of the $(F,E)$ adjunctions are weak equivalences, which we now do.

\begin{proposition}
If $\psi\colon\Oh\to\Cats$, then $\epsilon\colon FE\psi\to\psi$ is an objectwise
weak equivalence.
\end{proposition}

\begin{proof}
In fact, we shall show that, for every object $b$ of $\Oh$, the functor
$\epsilon\colon (FE\psi)(b) \to \psi(b)$ is a deformation retraction.
Let $\E = (FE\psi)(b)$. Tracing through the definitions, we see that
an object of $\E$ is a pair $(\beta,x)$ where $\beta\colon b\to c$ in $\Oh$
and $x\in\psi(c)$.
A map $(\beta,x)\to (\gamma\colon b\to d,y)$ is given by a pair
$(\delta,\zeta)$ where $\delta\colon c\to d$ with $\gamma = \delta\beta$
and $\zeta\colon x\to \delta^*y$.
With this notation, $\epsilon$ is given by
$\epsilon(\beta,x) = \beta^*x$ and $\epsilon(\delta,\zeta) = \beta^*\zeta$.

Now define $\xi\colon \psi(b)\to\E$ by
$\xi(z) = (1_b,z)$ on objects and
$\xi(\zeta) = (1_b,\zeta)$ on maps.
We have $\epsilon\xi = 1$, and we can define a
natural transformation $\omega\colon \xi\epsilon\to 1$ by
$\omega(\beta,x) = (\beta,1_{\beta^*x})$.

Thus, $\xi\colon\psi(b)\to (FE\psi)(b)$ is the inclusion of a deformation
retract and $\epsilon$ is a weak equivalence.
\end{proof}

\begin{proposition}\label{prop:counitequivalence}
If $\psi\colon \Oh\to\SSets$, then $\epsilon\colon FE\psi\to\psi$
is an objectwise weak equivalence.
\end{proposition}

\begin{proof}
This statement is part of \cite[1.4]{FG:diagrams}, but we reprove it here.
Let $b$ be an object of $\Oh$. We begin by describing the simplicial set
$(FE\psi)(b)$. We have the following pullback diagram.
\[
 \xymatrix{
  (FE\psi)(b) \ar[r] \ar[d] & \bar E\psi \ar[d] \\
  N(b\Over\Oh) \ar[r] & N\Oh
  }
\]
Thus, an $n$-simplex of $(FE\psi)(b)$ is given by two compatible maps
$\omega\colon\Delta[n]\to N(b\Over\Oh)$ and $\alpha\colon\Delta[n]\to \bar E\psi$.
The map $\omega$ is determined by a sequence
\[
 b\xrightarrow{g_0} \sigma_0\xrightarrow{g_1}\cdots\xrightarrow{g_n}\sigma_n
\]
where $\sigma_0\to\cdots\to\sigma_n$ is the composite map $\sigma\colon\Delta[n]\to N\Oh$.
The map $\alpha$ is then an element of $(E\psi)_\sigma$.
By Corollary~\ref{cor:fsigma} we can describe such an element by a sequence
$(\tau_0,\tau_1,\dots,\tau_n)$ where
$\tau_k\in \psi(\sigma_k)_{n-k}$ and
$d_0\tau_{k-1} = g_k^*\tau_k$ for $1\leq k \leq n$.
In summary, we can write an $n$-simplex of $(FE\psi)(b)$ as
\[
 (b\to\sigma_0\to\cdots\to\sigma_n; \tau_0,\dots,\tau_n)
\]
with the compatibility conditions above.

The map $\epsilon\colon (FE\psi)(b)\to \psi(b)$ is then given by
\[
 \epsilon(b\xrightarrow{g_0}\sigma_0\to\cdots\to\sigma_n; \tau_0,\dots,\tau_n)
  = g_0^*(\tau_0).
\]
We define $\xi\colon \psi(b)\to (FE\psi)(b)$ on an $n$-simplex $\tau$ by
\[
 \xi(\tau) = (b\to b\to\cdots\to b; \tau_0, d_0\tau_0,\dots, d_0^n\tau_0)
\]
where the maps $b\to b$ are all the identity.
Clearly $\epsilon\xi = 1$, and it is straightforward to define a simplicial
homotopy from $\xi\epsilon$ to the identity
(the construction is very similar to standard arguments involving the
bar construction).
Thus $\epsilon$ is a weak homotopy equivalence.
\end{proof}

From these two propositions and Proposition~\ref{prop:equivalence} we deduce
the following.

\begin{theorem}
With weak equivalence being weak fibre homotopy equivalence, 
$\Cats\Over\Oh$ has a homotopy category,
and $(F,E)\colon \Cats\Over\Oh\to \OCats$ induce inverse equivalences of homotopy
categories. Similarly, $\SSets\Over N\Oh$ has a homotopy category and
$(F,E)\colon \SSets\Over N\Oh\to \OSSets$ induce inverse equivalences of homotopy
categories.
\end{theorem}

Combining this theorem with Theorem~\ref{thm:diagramequivalence} we get
the following (cf. \cite[1.5]{FG:diagrams}).

\begin{corollary}
Each functor in the following commutative diagram induces an
equivalence of homotopy categories.
\[
 \xymatrix{
  \Cats\Over\Oh \ar[r]^F \ar[d]_N & \OCats \ar[d]^N \\
  \SSets\Over N\Oh \ar[r]_F & \OSSets
  }
\]
\end{corollary}

So far we have avoided discussing model category structures for
$\Cats\Over\Oh$ and $\SSets\Over N\Oh$. In the remainder of this paper
we construct a model category structure for $\SSets\Over N\Oh$
making $(F,E)\colon\SSets\Over N\Oh\to \OSSets$ a Quillen equivalence.
Given the difficult nature of Thomason's model category structure on
$\Cats$, and of any of the alternatives mentioned at the beginning of this
section, we shall not here attempt to construct a model category
structure on $\Cats\Over\Oh$.

\section{Right fibrations}

Before constructing the model category structure
on $\SSets\Over N\Oh$ we need to examine a weakening of the idea
of Kan fibration \cite{Kan:css} suggested by Proposition~\ref{prop:we}
and its failure in the case $k=0$.

\begin{definition}
Let $J_r = \{ \horn k n \to \Delta[n] \mid n\geq 1 \text{ and } 0 < k \leq n \}$.
A map of simplicial sets $p\colon X\to Y$ is a
{\em right fibration} if it has the right lifting property (RLP) with
respect to $J_r$ (i.e., is in $J_r$-inj, using the notation of
\cite{Hov:models}).
A map $i\colon A\to B$ is a {\em right anodyne extension}
if it has the left lifting property (LLP) with respect to all
right fibrations (i.e., is in $J_r$-cof).
\end{definition}

We could define analogous left fibrations and left anodyne extensions,
using $J_l = \{ \horn k n \to \Delta[n] \mid n\geq 1 \text{ and } 0 \leq k < n \}$.
However, we will have no use for these maps in this paper.
A (Kan) fibration
is a map that is both a right and a left fibration, so might be called
a two-sided fibration.

\begin{remark}
The use of the terms ``right'' and ``left'' here is, of course, somewhat arbitrary.
It is suggested by the
translation to categories that occurs on applying categorical realization
$c$ (the left adjoint to the nerve).
For categories, a functor $\phi\colon\C\to\D$
has the right lifting property with respect to
$c\horn 2 2\to [2]$ if the following is true: For every
pair of maps $\alpha\colon x\to z$ and $\beta\colon y\to z$ in $\C$ such that
there is a $g\colon \phi(x)\to\phi(y)$ such that $\phi(\alpha) = \phi(\beta)\circ g$,
there exists a $\gamma\colon x\to y$ such that $\phi(\gamma) = g$
and $\alpha = \beta\circ\gamma$. We call this condition
``right divisibility'' (the right factor $\gamma$ can be found).
Similarly, having the right lifting property with respect to
$c\horn 0 2\to [2]$ we call left divisibility.
A right fibration, therefore, has right divisibility but not necessarily
left divisibility.
\end{remark}

The results of \cite[\S 3.3]{Hov:models} are true with anodyne extensions
replaced by right anodyne extensions and fibrations replaced by
right fibrations; Hovey's proofs are easily modified for this case.
In particular, the two main results of that section become the following.

\begin{theorem}\label{thm:boxprod}
If $i\colon K\to L$ is an inclusion of simplicial sets
and $f\colon A\to B$ is a right anodyne extension, then
\begin{equation*}
 i\boxprod f \colon (K\times B)\disjunion_{K\times A}(L\times A) \to L\times B
\end{equation*}
is a right anodyne extension.
\end{theorem}

\begin{theorem}\label{thm:mapbox}
If $i\colon K\to L$ is an inclusion of simplicial sets
and $p\colon X\to Y$ is a right fibration, then
\begin{equation*}
 \Map_\boxprod(i,p)\colon \Map(L,X)\to \Map(K,X)\times_{\Map(K,Y)}\Map(L,Y)
\end{equation*}
is a right fibration.
\end{theorem}

\begin{remark}
It's interesting that we can carry through the definition of homotopy groups
as in \cite[\S 3.4]{Hov:models} for right (or left) fibrant simplicial sets
as easily as for fibrant ones (we just have to modify the proof that homotopy of
vertices is an equivalence relation to avoid using the missing horns).
By comparing the homotopy groups of $X$ with those of $\Omega X$
(as in \cite[\S 3.6]{Hov:models}) we can show that the homotopy groups we
define this way agree with the usual ones
(the homotopy groups of the geometric realization).
However, we shall make no direct use of homotopy groups in this paper.
\end{remark}

There is another construction that gives us right anodyne extensions.
To describe it, we first discuss simplicial cones.
Let $\Delta^-$ denote the category whose objects are the objects
of $\Delta$ together with one extra object, the empty set.
For notational convenience we write $[-1] = \emptyset$.
Maps are again order-preserving functions, so
$[-1]$ is initial in $\Delta^-$. We can extend any
simplicial set $A$ naturally to a functor on $\Delta^-$ by
letting $A_{-1} = *$, a one-point set.
$\Delta[-1]$ is the empty simplicial set.
($\Delta^-$ is the category called $\Delta$ in
\cite{Mac:categories}; contravariant set-valued functors on
$\Delta^-$ are there called augmented simplicial sets.)
Let $\SSets_*$ denote the category of based simplicial sets.

\begin{definition}
Let $C\Delta\colon \Delta^-\to \SSets_*$ be the covariant functor with
$C\Delta[n] = \Delta[n+1]$ with $0$ as the base vertex. 
If $f\colon[m]\to[n]$,
let $\bar f\colon [m+1]\to[n+1]$ be the function defined by
$\bar f(0) = 0$ and $\bar f(k) = f(k-1) + 1$ for $k>0$;
let $C\Delta(f) = \bar f^*$.
Now extend $C$ to a functor $C\colon\SSets\to\SSets_*$ by letting
\[
 CA = \int^{n\in\Delta^-} \!\!\!\!A_n \times C\Delta[n].
\]
\end{definition}

The effect of $C$ is to add an initial ``cone point'' to every simplex
of $A$, these cone points all being identified to become the base vertex $*$
of $CA$.
There is a natural inclusion $j\colon A\to CA$ (of unbased simplicial sets),
induced by the map $\delta_0\colon \Delta[n]\to\Delta[n+1]$,
the inclusion of the face opposite $0$.

The description of the functor $C$ as a coend shows that it
has as its right adjoint $P\colon\SSets_*\to\SSets$, where
\[
 (PB)_n = \Hom_{\SSets_*}(C\Delta[n],B)
 = \{ z\colon\Delta[n+1]\to B \mid z(0) = * \}.
\]
Dual to the inclusion $j$ is the
(unbased) map
$\bar\jmath\colon PB\to B$ that takes an $(n+1)$-simplex in $B$ to its
face opposite 0.
Although we've described $PB$ as an unbased simplicial set,
it has a canonical base vertex $*$, given by the map
$C\Delta[0]\to B$ constant at the base vertex of $B$.

\begin{theorem}\label{thm:cones}
If $f\colon A\to B$ is an anodyne extension then
$\hat C f\colon B\disjunion_A CA \to CB$
is a right anodyne extension.
\end{theorem}

\begin{proof}
We first check the special case $f\colon\horn kn\to\Delta[n]$.
In this case, it is easy to see that
$\hat C f$ is the inclusion $\horn{k+1}{n+1}\to\Delta[n+1]$, which
is a right anodyne extension since $k+1>0$.
If we let
\[
 J = \{ \horn k n \to \Delta[n] \mid n\geq 1 \text{ and } 0 \leq k \leq n \},
\]
and let $\hat C J$ denote the set that results when we apply $\hat C$
to each map in $J$, then we have shown that
every map in $\hat C J$ is a $J_r$-cofibration.

To get to the general case we note first that, by adjunction,
liftings in a diagram as on the left below are in one-to-one
correspondence with liftings in the diagram on the right.
\[
 \xymatrix{
  B\disjunion_A CA \ar[d]_{\hat C f} \ar[r] & X \ar[d]^p \\
  CB \ar@{-->}[ur] \ar[r] & Y
  }
\qquad\qquad
 \xymatrix{
  A \ar[d]_f \ar[r] & PX \ar[d]^{\hat P p} \\
  B \ar@{-->}[ur] \ar[r] & PY\times_Y X
  }
\]
In particular, if $p$ has the RLP with respect to $\hat C J$,
then $\hat P p$ has the RLP with respect to $J$
(so is a Kan fibration).

Since every map in $\hat C J$ is a $J_r$-cofibration,
every map in $\hat C J$ has the LLP with respect to $J_r$-inj,
the set of $J_r$-injectives (the right fibrations).
By adjunction, every map in $J$ has the LLP with respect to
$\hat P (J_r\text{-inj})$.
It follows that every $J$-cofibration has the LLP with respect
to $\hat P (J_r\text{-inj})$, hence, by adjunction, that
every map in $\hat C(J\text{-cof})$ has the LLP with respect
to $J_r$-inj. But this says that, if $f$ is an anodyne
extension, then $\hat C J$ is a right anodyne extension.
\end{proof}

Cones are naturally contractible, in the following sense.

\begin{lemma}\label{lem:contraction}
There is a natural map $H\colon C(A\times\Delta[1])\to CA$
with the property that $H|C(A\times 0)$ is constant at the base vertex
while $H|C(A\times 1)$ is the identity.
\end{lemma}

\begin{proof}
We first define $H\colon C(\Delta[n]\times\Delta[1])\to C\Delta[n]$.
Recall that $\Delta[n]\times\Delta[1]$ has $n+1$ top-dimensional
nondegenerate simplices $v_0$ through $v_n$, where
$v_k$ is given by the sequences of vertices
$(0,0)$, $(1,0)$, \dots, $(k,0)$, $(k,1)$, \dots, $(n,1)$.
Corresponding to these are $n+1$ top-dimensional simplices of
$C(\Delta[n]\times\Delta[1])$ we shall call $w_0$ through $w_n$.
We define $H|w_k\colon\Delta[n+2]\to\Delta[n+1]$
to be $\sigma_0^{n+1-k}$.
We leave it to the reader to verify that this does define a simplicial
map $H$ with the stated properties, and that $H$ is natural on
$\Delta^-$.
We then extend $H$ to all of $\SSets$ using the coend:
\[
 C(A\times\Delta[1])
  \iso \int^n \!\!A_n\times C(\Delta[n]\times\Delta[1])
   \to \int^n \!\!A_n\times C\Delta[n] = CA.
\]
\end{proof}

By adjunction, we get a contraction of the dual construction.

\begin{corollary}\label{cor:Pcontraction}
There is a natural map $\hat H\colon PB\times\Delta[1]\to PB$
such that $\hat H|PB\times 0$ is constant at the base vertex
and $\hat H|PB\times 1$ is the identity.
\end{corollary}

A key technical result is the following.

\begin{lemma}\label{lem:acyclic}
Suppose that $p\colon X\to Y$ is a right fibration of simplicial sets.
Suppose also that, for every vertex $v\in Y$, the fiber $p^{-1}(v)$ is
an acyclic Kan complex. Then $p$ is an acyclic fibration.
\end{lemma}

\begin{proof}
We show that $p$ has the RLP with respect to every map
$\bndry\Delta[n]\to\Delta[n]$. Consider the following lifting
problem:
\[
 \xymatrix{
  \bndry\Delta[n] \ar[r]^-x \ar[d] & X \ar[d]^p \\
  \Delta[n] \ar@{-->}[ur] \ar[r]_-y & Y
  }
\]
That this lifting problem can be solved for $n=0$ is just the statement
that the fibers, being acyclic, are nonempty.
Let us now assume that $n>1$ and think of $\bndry\Delta[n]\to\Delta[n]$
as $\Delta[n-1]\union C\bndry\Delta[n-1]\to C\Delta[n-1]$. 
We shall solve the lifting problem
by first pulling it into a fiber, where the problem is solvable
by assumption.

To do this, we consider first the following lifting problem:
\[
 \xymatrix{
  (\Delta[n-1]\times 1)\union C(\bndry\Delta[n-1]\times 1) \ar[r]^-x \ar[d]
   & X \ar[d]^p \\
  (\Delta[n-1]\times\Delta[1])\union C(\bndry\Delta[n-1]\times\Delta[1])
   \ar[r]_-{yH} \ar@{-->}[ur]^h & Y
  }
\]
The map $yH$ is the restriction of the composite of $y$ with the contraction
$H\colon C(\Delta[n-1]\times\Delta[1])\to C\Delta[n-1]$
of Lemma~\ref{lem:contraction}.
We claim that the map on the left is a right anodyne extension. We see this
by writing it as the composite
\begin{align*}
 (\Delta[n-1]\times 1)\union C(\bndry\Delta[n-1]\times 1)
 &\to (\Delta[n-1]\times \Delta[1])\union C(\bndry\Delta[n-1]\times 1) \\
 &\to (\Delta[n-1]\times \Delta[1])\union C(\bndry\Delta[n-1]\times \Delta[1]).
\end{align*}
The first of these maps is obtained by pushing out along
$\Delta[n-1]\times 1 \to \Delta[n-1]\times\Delta[1]$, which is
a right anodyne extension by Theorem~\ref{thm:boxprod}.
The second map is obtained by pushing out along
\[
 (\bndry\Delta[n-1]\times\Delta[1])\union C(\bndry\Delta[n-1]\times 1)
  \to C(\bndry\Delta[n-1]\times\Delta[1]),
\]
which is a right anodyne extension by Theorem~\ref{thm:cones}.
Therefore, the composite is a right anodyne extension as claimed, and
so we can find the map $h$ in the diagram above.

Now note that the image of 
$h|(\Delta[n-1]\times 0)\union C(\bndry\Delta[n-1]\times 0)$
lies entirely in a fiber
of $p$. By assumption, then, we can extend $h$ over
$C(\Delta[n-1]\times 0)$. We now consider the following lifting problem:
\[
 \xymatrix{
  (\Delta[n-1]\times\Delta[1])\union
   C(\Delta[n-1]\times 0\union\bndry\Delta[n-1]\times\Delta[1])
    \ar[r]^-h \ar[d] & X \ar[d]^p \\
  C(\Delta[n-1]\times\Delta[1]) \ar[r]_-{yH} \ar@{-->}[ur]^k & Y
  }
\]
The map on the left is a right anodyne extension by
Theorem~\ref{thm:cones}, so we can find the lift $k$.
The restriction of $k$ to $C(\Delta[n-1]\times 1)$ is then the solution
to our original lifting problem.
\end{proof}

We also need to examine the relation between fibers and homotopy fibers.

\begin{definition}
If $Y$ is a simplicial set and $y\in Y_0$ is a vertex, let
$(Y,y)$ denote the based simplicial set we get by considering
$y$ as the base vertex of $Y$. $P(Y,y)$ is then 
the simplicial set with
\begin{equation*}
 P(Y,y)_n = \{ z\colon \Delta[n+1]\to Y \mid z(0) = y \}.
\end{equation*}
If $p\colon X\to Y$ is a map of simplicial sets and $y\in Y_0$
is a vertex, we define $P(p,y)$,
the {\em (right) homotopy fiber over $y$}, to be the pullback in the
following diagram:
\[
 \xymatrix{
  P(p,y) \ar[r] \ar[d]_{\bar p} & X \ar[d]^p \\
  P(Y,y) \ar[r]_{\bar\jmath} & Y
  }
\]
\end{definition}

When $Y = N\Oh$, $P(Y,y) = N(y\Over\Oh)$ and so
$P(p,y) = (Fp)(y)$, as defined in Section~\ref{sec:functors}.

Recall that the fiber $p^{-1}(y)$ is the pullback in the following diagram:
\[
 \xymatrix{
  p^{-1}(y) \ar[r] \ar[d] & X \ar[d]^p \\
  \Delta[0] \ar[r]_y & Y
  }
\]
The map $*\colon \Delta[0]\to P(Y,y)$ then induces a map
$p^{-1}(y) \to P(p,y)$.

\begin{proposition}\label{prop:fibers}
If $p\colon X\to Y$ is a right fibration of simplicial sets and
$y$ is a vertex of $Y$, then the map
$p^{-1}(y)\to P(p,y)$ is the inclusion of a deformation retract.
\end{proposition}

\begin{proof}
Since right lifting properties are preserved by pullback and $p\colon X\to Y$
is a right fibration, so is $\bar p\colon P(p,y) \to P(Y,y)$.
Let $\hat H\colon P(Y,y)\times\Delta[1]\to P(Y,y)$ be the deformation
constructed in Corollary~\ref{cor:Pcontraction} and consider the following lifting
problem.
\[
 \xymatrix{
  p^{-1}(y)\times\Delta[1] \union P(p,y)\times 1 \ar[r]\ar[d]
   & P(p,y) \ar[d]^{\bar p} \\
  P(p,y)\times\Delta[1] \ar@{-->}[ur] \ar[r]_-{\hat H\circ \bar p}
   & P(Y,y)
  }
\]
Since $p^{-1}(y)\to P(p,y)$ is an inclusion, it follows from
Theorem~\ref{thm:boxprod} that
the map on the left
is a right anodyne extension, hence we can solve this lifting problem.
The lifted map is the desired deformation.
\end{proof}

\begin{remark}
This result is clearly related to Proposition~\ref{prop:counitequivalence}.
However, it is not clear what set of assumptions would imply both results.
Neither is it clear that such a unification would be useful here.
\end{remark}

\section{The model category structure and the Quillen equivalence}\label{sec:model}

We now define a new model category structure on $\SSets\Over N\Oh$.

\begin{definition}
Let $\phi\colon X\to N\Oh$ and $\psi\colon Y\to N\Oh$ be simplicial sets
over $N\Oh$, and let $f\colon X\to Y$ be a map over $N\Oh$.
\begin{enumerate}
\item Say that $f$ is a {\em weak equivalence over $N\Oh$}
 if, for each vertex $b\in N\Oh$,
 $(F\phi)(b)\to (F\psi)(b)$ is a weak equivalence of simplicial sets.
\item Say that $f$ is a {\em cofibration over $N\Oh$} if it is a cofibration
 of simplicial sets, i.e., an inclusion.
\item Say that $f$ is a {\em fibration over $N\Oh$} if it is a right fibration
 and if, for each vertex $b\in N\Oh$, $\phi^{-1}(b)\to \psi^{-1}(b)$ is a fibration of
 simplicial sets.
\end{enumerate}
\end{definition}

\begin{theorem}\label{thm:modelstructure}
With the definitions of weak equivalence, cofibration, and fibration over $N\Oh$ given
above, $\SSets\Over N\Oh$ is a cofibrantly generated model category. For a set of
generating cofibrations we may take the set of all inclusions
$\bndry\Delta[n]\to\Delta[n]$ over $N\Oh$.
For a set of generating acyclic cofibrations we may take the set of
all inclusions $\horn k n\to\Delta[n]$ over $N\Oh$ with $k>0$, together with the
inclusions $\horn 0 n\to\Delta[n]$ in which the map $\Delta[n]\to N\Oh$ is constant
at a vertex.
\end{theorem}

\begin{proof}
We shall use the ``recognition principle'' \cite[2.1.19]{Hov:models}
for cofibrantly generated categories. Hovey gives six criteria to check to
verify that there is a cofibrantly generated model category structure with the
stated weak equivalences, with the stated set $I$ as a set of generating cofibrations and
the stated set $J$ as a set of generating acyclic cofibrations.

(1) {\em The collection of weak equivalences has the two-out-of-three property
and is closed under retracts.} This is straightforward to check given that it is true for
ordinary weak equivalences of simplicial sets.

(2) \& (3) {\em The domains of $I$ and $J$ are small.} This follows from the fact that
all simplicial sets are small.

(4) {\em Relative $J$-cell complexes are acyclic $I$-cofibrations.}
Every map in $J$ is an inclusion, from which it follows that relative $J$-cell complexes
are $I$-cofibrations. It suffices, then, to check that every map in $J$ is a
weak equivalence over $N\Oh$. But this is the content of Proposition~\ref{prop:we}.
(The case $\horn 0 n \to\Delta[n]\to\text{vertex}$ is easy.)

(5) {\em $I$-injectives are weak equivalences.}
Let $\phi\colon X\to N\Oh$ and $\psi\colon Y\to N\Oh$, and
let $p\colon X\to Y$ be an $I$-injective over $N\Oh$. Then $p$ is an
acyclic Kan fibration. If $b$ is any vertex in $N\Oh$, a diagram chase shows that
the following is a pullback diagram.
\[
 \xymatrix{
  (F\phi)(b) \ar[r] \ar[d] & X \ar[d]^p \\
  (F\psi)(b) \ar[r] & Y
  }
\]
It follows that $(F\phi)(b)\to (F\psi)(b)$ is an acyclic Kan fibration
for every $b$, hence that
$X\to Y$ is a weak equivalence over $N\Oh$.

(6) {\em Acyclic $J$-injectives are $I$-injectives.}
Let $\phi\colon X\to N\Oh$ and $\psi\colon Y\to N\Oh$, and
let $p\colon X\to Y$ be an acyclic $J$-injective over $N\Oh$.
By the inclusion in $J$ of the maps $\horn 0 n\to \Delta[n]\to\text{vertex}$,
it follows that $\phi^{-1}(b)\to \psi^{-1}(b)$ is a Kan fibration for every
vertex $b$ in $N\Oh$. On the other hand, we have the following commutative
diagram, in which the horizontal maps are inclusions of deformation
retracts (by Proposition~\ref{prop:fibers})
and the map on the right is a weak equivalence by assumption.
\[
 \xymatrix{
  \phi^{-1}(b) \ar[r] \ar[d] & (F\phi)(b) \ar[d] \\
  \psi^{-1}(b) \ar[r] & (F\psi)(b)
  }
\]
It follows that $\phi^{-1}(b)\to \psi^{-1}(b)$ is a weak equivalence, hence
an acyclic Kan fibration, for each $b$.
From this we can conclude that, for each vertex $y$ of $Y$,
$p^{-1}(y)$ is an acyclic Kan complex.
We can now appeal to Lemma~\ref{lem:acyclic} to conclude that
$p$ is an acyclic Kan fibration, hence is an $I$-injective.

Thus, by \cite[2.1.19]{Hov:models}, $\SSets\Over N\Oh$ is a cofibrantly
generated model category with the stated weak equivalences and generating cofibrations.
What remains to show is the characterizations of cofibrations and fibrations
over $N\Oh$. But, it is obvious that the $I$-cofibrations are just the usual
cofibrations of simplicial sets, i.e., inclusions. It is also clear that
the $J$-injectives are precisely the right fibrations for which each
map of fibers $\phi^{-1}(b)\to \psi^{-1}(b)$ is a Kan fibration.
\end{proof}

Note that we can also characterize acyclic fibrations over
$N\Oh$: They are precisely the usual acyclic Kan fibrations because this
is the class that has the RLP with respect to $I$.

Now we want to compare this model category to $\OSSets$. 
As mentioned earlier,
the model category structure on $\OSSets$ is cofibrantly generated;
for a set of generating cofibrations we may take the collection of maps
$\Oh(-,c)\times\bndry\Delta[n]\to\Oh(-,c)\times\Delta[n]$,
and for a set of generating acyclic cofibrations we may take the collection
of maps
$\Oh(-,c)\times\horn kn \to\Oh(-,c)\times\Delta[n]$, for all $k$ and $n$.

\begin{theorem}\label{thm:equivalence}
Let $\SSets\Over N\Oh$ and $\OSSets$ have the model category structures
described above. Then $(F,E)\colon \SSets\Over N\Oh \to \OSSets$
is a Quillen equivalence.
\end{theorem}

\begin{proof}
We show first that $F$ preserves cofibrations and acyclic cofibrations.
To show that $F$ preserves cofibrations, it suffices to show that it takes
each generating cofibration to a cofibration. But, Lemma~\ref{lem:attachment}
shows that, if $\sigma\colon\Delta[n]\to N\Oh$ and $\lambda$ is the
composite $\bndry\Delta[n]\to\Delta[n]\to N\Oh$, we
have the following pushout diagram.
\[
 \xymatrix{
  \Oh(-,\sigma_0)\times\bndry\Delta[n] \ar[r] \ar[d] & F\lambda \ar[d]\\
  \Oh(-,\sigma_0)\times\Delta[n] \ar[r] & F\sigma
  }
\]
From this it is clear that $F\lambda\to F\sigma$ is a cofibration,
so $F$ preserves cofibrations.
On the other hand, $F$ preserves weak equivalences by definition,
so $F$ also preserves acyclic cofibrations.
Therefore, $(F,E)$ is a Quillen adjunction.

To show that $(F,E)$ is a Quillen equivalence, we need to show that,
if $\phi\colon A\to N\Oh$ is cofibrant (which is always true) and
$\psi\colon \Oh\to\SSets$ is fibrant, then a map
$f\colon \phi \to E\psi$ is a weak equivalence if and only if its adjoint
$\hat f\colon F\phi\to \psi$ is.
Since $f$ is a weak equivalence if and only if $Ff$ is,
it suffices to show that the counit $\epsilon\colon FE\psi\to \psi$
is a weak equivalence, but this was done in Proposition~\ref{prop:counitequivalence}.
\end{proof}

In fact, $\SSets\Over N\Oh$ and $\OSSets$ are both simplicial model categories
and $(F,E)$ is a Quillen equivalence of simplicial model categories, as can
be easily checked.

\providecommand{\bysame}{\leavevmode\hbox to3em{\hrulefill}\thinspace}

\end{document}